\documentclass[12pt]{article}
\usepackage[letterpaper,margin=1in]{geometry}
\usepackage[T1]{fontenc}
\IfFileExists{newtxtext.sty}{\usepackage{newtxtext,newtxmath}}{\usepackage{mathptmx}}
\usepackage{setspace}
\usepackage{microtype}
\usepackage{url}
\newcommand{\fasection}[1]{%
\vspace{1.4em}
{\centering\textbf{#1}\par}
\vspace{0.6em}\noindent}

\begin{document}

\begin{center}
{\LARGE\bfseries The Wandering Must Continue}\\[0.8em]
{\large\itshape Pure Mathematicians as Rangers of the North}\\[1.2em]
{\normalsize\scshape Jun--Yong Park}
\end{center}

\fasection{I. THE SEMINAR LIST}

I have a test. When I travel to a country I have never visited and want to understand its scientific and technological capacity --- not its GDP, not its venture capital, not its startup culture --- I go to the best university in that country. I walk into the mathematics department. I find the seminar schedule, usually posted on a board near the entrance or taped to a wall on the first floor. And I read the titles.

\smallskip

That is all I need.

\smallskip

If the titles read something like ``Distribution of Hodge loci and atypical special points'' or ``Resolution of singularities in positive characteristic'' or ``Unboundedness of Mordell--Weil rank over global fields,'' I know that this country maintains deep scientific infrastructure. Not because anyone will ever apply these results to build a bridge or cure a disease --- almost certainly no one will --- but because the presence of people doing this kind of work is a reliable signal that the broader ecosystem is healthy. Good physicists were trained here. Good engineers passed through the same corridors. Good AI researchers working on alignment took the same foundational courses.

\smallskip

Applied titles tell me other things, also worth knowing: how deeply the country's scientific expertise reaches, what level of analysis its ministries bring to strategy, where the ceiling lies on nonlinear physical modeling built on rigorous numerics. Between pure and applied mathematics runs a relay, not a wall, passing ideas in both directions. But it is the pure titles that make the test work, for a simple reason: nobody has an instrumental motive to fake them. No ministry ever demanded a seminar on Hodge loci. No investor was ever impressed by one. If it is on the board, it is there because the culture that produces such things is alive in the building. That is what makes the signal honest.

\smallskip

This is an observation about national infrastructure. The seminar list is a barometer, and like all barometers, it measures something other than what it displays. It displays abstruse mathematics. It measures the depth of a civilization's capacity to reason about the structures that underlie its own science and technology.

\fasection{II. RANGERS OF THE NORTH}

Tolkien understood something about invisible protection.

\smallskip

In \textit{The Lord of the Rings}~\cite{tolkien1954}, the Shire, a green, peaceful corner of Middle-earth, is kept safe not by its own inhabitants but by the Rangers: grim, weathered men who patrol its borders in silence, fighting threats that the hobbits have never seen and do not believe exist. The hobbits regard the Rangers with suspicion and faint contempt. ``He is one of the wandering folk --- Rangers we call them,'' says the innkeeper at Bree, in a tone that suggests vagrancy more than vigilance.

\smallskip

The Rangers do not glitter; they produce nothing visible. They wander the wild lands in the dark so that inside the Shire, ordinary life can proceed as if order were the natural state of the world. Because every threat is met miles beyond the hedgerows, the very success of the patrol erases the evidence that a patrol was ever needed.

\smallskip

I think about the Rangers more often than I should, because they describe, with uncomfortable precision, the institutional position of pure mathematicians. When a civilization's conceptual borders have been secure for generations, it is easy to believe that rigor maintains itself, that definitions stay sharp on their own, and that the people wandering the abstract frontier are superfluous.

\fasection{III. WHAT MATHEMATICIANS ACTUALLY DO}

The public image of a mathematician is a person who solves very hard problems. This is not entirely wrong, but it is wrong in the way that matters. It is like saying that a physician is a person who writes prescriptions: technically true, and fundamentally misleading.

\smallskip

What mathematicians actually do, most of the time, is maintain a culture. They sustain a set of practices, habits of thought, standards of argument, and norms of honesty that have been built up over centuries and that cannot be written down in any manual. A graduate student does not simply learn how to prove new theorems. She learns how to be suspicious of her own reasoning. She learns that an argument that ``feels right'' is not right until every step has been verified, and that stating a result precisely is as important as proving it. She learns to be grateful when someone finds an error in her work, because the alternative, an error that persists, is worse.

\smallskip

But that discipline is only the preparation for a deeper transformation. At some point in her apprenticeship, a threshold is crossed: she stops experiencing mathematics as a game of rules played with symbols on paper and begins to see the structures themselves. She learns to reach for the right geometric principle or algebraic structure, so that a new theorem does not merely stand as an isolated fact, but illuminates the landscape around it. What once looked like formal manipulation becomes an act of perception, an encounter with objects that have their own stubborn necessity. She learns that the difference between a good mathematician and a mediocre one is not speed or cleverness but a quality of attention and depth of awareness: a refusal to let anything pass that has not been fully understood, or to compromise true understanding for a mechanical answer.

\smallskip

These habits of mind and sight are not natural. They are not inevitable. They must be cultivated, by example, over years, in academic communities of people who share them.

\smallskip

And these habits are, I submit, among the most valuable cultural artifacts that human civilization has produced. Their value lies not in the theorems they generate, most of which are read by a few dozen people and forgotten, but in the habits of rigor, courage, and precision that radiate outward from mathematics departments into every field they touch. You see their trace in the physicist who distrusts a convenient heuristic until the derivation is sound, the engineer who refuses to shave a structural safety margin, and the AI researcher who checks every edge case rather than assuming the pattern holds. These instincts were shaped, directly or indirectly, by contact with a culture that treats the difference between ``almost certainly true'' and ``proved'' as the difference between nothing and everything.

\smallskip

The patrol that mathematicians conduct is invisible to everyone outside the community. Its work is quiet, and it has no market value or political constituency. Like clean water or functioning courts, it is the kind of invisible infrastructure a society notices only after it is gone and the wellspring of real discovery runs dry.

\fasection{IV. WHAT CANNOT BE RECONSTITUTED}

The danger is not that mathematicians will be replaced overnight; it is slower, quieter, and for that reason harder to resist.

\smallskip

A university facing budget pressure does not abolish its mathematics department. It freezes a position, drops real and complex analysis from the undergraduate core, and eliminates one-year foundational courses in algebraic geometry and topology. The seminar list gets a little shorter. The graduate program admits one fewer student. A postdoctoral fellowship is converted to a teaching position. Each individual decision is rational. Each is, on its own, minor.

\smallskip

The result, compounded over a decade, is that this culture begins to thin. Not to vanish entirely, but to thin, the way a forest thins when its oldest trees are not replaced. You can still walk through it. You can still call it a forest. But its capacity to sustain the ecosystem it once supported has quietly diminished.

\smallskip

This process is already underway. We see it in quiet, localized attrition and, increasingly, in structural amputation. In September 2023, West Virginia University, a flagship state research institution, voted to discontinue its master's and doctoral mathematics programs to close a budget deficit~\cite{wvu2023}. The mechanism is not dramatic: no one is burning libraries. The mechanism is arithmetic: fewer positions, fewer students, fewer grants, fewer seminars. The barometer drops.

\smallskip

The critical fact about mathematical capacity is that it cannot be reconstituted on demand. You cannot mass produce a mathematician through an accelerated three-month training program, or conjure a research culture the way you scale a data center. The training takes ten to fifteen years and depends on a tradition of apprenticeship that stretches back centuries, because what is being transmitted is not a manual of procedures, but a way of seeing, passed from one mind to another.

\smallskip

It took three thousand years to build this chain. It would take a single generation to break it, and perhaps a century to notice what was lost.

\fasection{V. THE CERTIFICATE AND THE CANDLE}

Here I must answer the obvious objection, because every technically literate reader has already raised it. We are entering an era where artificial intelligence does not merely verify mathematics, but autonomously discovers it. OpenAI has recently announced an AI-generated proof that the three-dimensional incompressible Navier--Stokes equations, with smooth external forcing, can develop a singularity in finite time, and has released a formalization in Lean~\cite{openai2026}. Formal proof assistants check every step of an argument against the axioms, mechanically and without fatigue. If AI systems produce the theorems and proof assistants certify them, the pipeline from conjecture to formally verified result appears complete. Why maintain, at enormous expense, a guild of artisanal mathematicians?

\smallskip

It is a good objection, and I concede everything in it except the conclusion.

\smallskip

A proof certificate establishes that a conclusion follows from the formal assumptions. By itself, it need not make the argument's central ideas intelligible, because the system that produced it requires no inner conceptual horizon in which mathematical objects are inhabited and seen as a unified landscape. A calculation can verify by direct addition that the numbers from $1$ to $100$ sum to $5050$; understanding begins only when one writes the sequence forward and backward, pairs the columns into a hundred copies of $101$, and sees the symmetry that makes the sum inevitable. When an autonomous system produces a large, independently checked formal proof establishing a finite-time singularity under specified hypotheses, we can have every reason to trust the certificate while still lacking an illuminating account of the mechanism or its broader implications. What threatens the discipline is not that neural networks can generate proofs and formal systems can check them, but the economic conclusion that administrators and funding bodies will be tempted to draw: that because generation and verification can be automated, human comprehension no longer needs to be paid for. Establishing a result and understanding why it holds are distinct achievements.

\smallskip

There is a deeper problem at the boundary between formalism and reality: a formal system can verify syntax, but it cannot verify meaning. A proof assistant will happily certify a flawless proof of a subtly different theorem. Deciding what to prove, and recognizing that a formal statement means what we intend it to mean --- that the tangible string of symbols captures the mathematical reality at stake --- is not clerical work, but an exercise of human judgment. It is often the hardest part of the subject, and no certificate can contain it. That judgment is what connects the formalism to the world.

\smallskip

And it is at the connection to the world that things fail. When a bridge collapses, when a global financial system implodes, when an artificial superintelligence is misaligned, the question is never ``was the certificate valid?'' The certificate was valid. The question is: what did the model leave out? What was assumed away, and what was optimized at the expense of what? No system can be the sole auditor of its own presuppositions. Someone standing outside it must be able to ask the right question: a human being carrying with her three thousand years of tradition, with the training and the discernment to find, in a chain of impeccable reasoning, the one link that was attached to nothing.

\smallskip

The astronomer Carl Sagan called scientific understanding a candle in the dark~\cite{sagan1995}, our protection against a world of technology we use but cannot evaluate, governed by claims we cannot test. We are now building that world by design. We are rapidly approaching a condition where we will live surrounded by marvelous devices whose operation we do not comprehend, based on engineering principles we cannot audit, derived from formalisms we cannot interpret. When everything works, this is merely eerie. When something fails, it is the difference between a diagnosis and a mystery.

\smallskip

This is the service that mathematical culture provides: not merely theorems, but a standard of conceptual clarity. It is the trained capacity to assess mathematical claims, examine their assumptions, and explain what they do and do not establish; to distinguish the rigorous from the plausible; and to insist, stubbornly, that ``it works'' is not the same as ``we understand why it works.''

\smallskip

Beyond verification, large language models excel at searching within a settled formal system once the axioms and the target conjecture have already been laid down. What they cannot do on their own is recognize when an unformalized mystery, whether in the architecture of neural populations or the latent representation spaces of frontier AI systems, cannot be captured by nineteenth-century calculus and linear probes, and requires the \textit{structural grammar} of modern mathematics. Today's tools reveal fragments: circuits, features, scaling curves. But fragments are not a theory. We can inspect portions of the machine without knowing what sort of mathematical object the whole machine has become. A machine can brute-force a puzzle inside an existing box; it takes a human mind immersed in foundational structure, and the inner horizon it cultivates, to see which field of pure mathematics is the true native language of a new domain of science.

\fasection{VI. THE PATROL MUST CONTINUE}

Pure mathematics looks, from the outside, like a spectacular waste of talent and resources: thousands of highly trained minds engaged in activities of zero practical consequence. Its products are arcane. Its practitioners spend years on problems whose solutions will be celebrated by a handful of specialists and unrecognized by everyone else.

\smallskip

And this is precisely why it works.

\smallskip

The apparent ``uselessness'' of pure mathematics is the condition of its strength. It is what ensures that rigor, patience, and the subordination of wishful thinking to proof are cultivated for their own sake, free from immediate commercial incentives. A mathematician who proves a theorem about the distribution of prime numbers is not trying to break anyone's encryption. She is trying to see the truth of arithmetic because her insights led her there. The fact that her understanding might, decades later, be relevant to encryption is a consequence of the depth of her inquiry, not its purpose.

\smallskip

If you redirect this community toward ``useful'' problems, if you fund only research with identifiable economic impact, if you replace exploratory inquiry with targeted problem-solving, you do not get more useful mathematics. You get shallower mathematics. You get mathematics that is optimized for the next grant cycle rather than the next century. You get mathematics that has lost the very quality, its disinterested depth, that made it useful in the first place.

\smallskip

And you get a generation of students who have learned to solve problems, but not to think; who can operate a system, but not see its limitations; who can follow rules, but not imagine alternatives.

\smallskip

The Rangers' patrol, the deliberate wandering, must continue, not because any single path through the wild matters on its own, but because the patrol itself is what keeps the borders safe.

\fasection{VII. WHAT IS TO BE DONE}

The world we live in changes through capital, policy, and law. It is institutions that sustain the discipline, and any serious defense of it must reckon with the realities facing the elected officials who oversee government funding, the university administrators tasked with balancing budgets, and above all the taxpaying public, who trust our standards. Elsewhere, I have proposed specific mechanisms for preserving mathematical capacity~\cite{park2026}; expanded here, four of them form the core of what is needed now.

\smallskip

First, \textbf{independent offices of mathematical assurance}: standing public bodies, funded by governments rather than by AI companies, whose mathematicians examine the foundational logic and geometric architecture of machine-generated reasoning wherever it enters consequential systems, including infrastructure, cryptography, medicine, and finance. Mathematicians are not there to replace physicians or civil engineers, but to audit the structural mathematics upon which those specialists increasingly rely. The model is the air-safety investigator: independent, uncompromising, and answerable to the public.

\smallskip

Second, \textbf{comprehension residencies at mathematical institutes}: funded positions devoted to understanding major machine-produced proofs, reconstructing them, simplifying them, and extracting the ideas that a bare certificate need not make intelligible. Here the deliverable is not a publication, but an explanation.

\smallskip

Third, \textbf{cross-disciplinary structural fellowships}: joint initiatives that embed pure mathematicians directly alongside neuroscientists, physicians, physicists, and engineers, using modern computational and AI tools to translate deep geometric and algebraic frameworks into predictive models of complex systems. Vast territories of twenty-first-century biology, medicine, and artificial intelligence remain bottlenecked not by a lack of data, but by a lack of mathematical architecture; the pure mathematician brings the structural language that empirical science cannot extract from data alone.

\smallskip

Fourth, \textbf{reinforced graduate programs}: training that fiercely protects year-long foundational courses in pure mathematics while sustaining a seminar culture at every level, from student reading groups to invited colloquia and faculty-organized research seminars. The ambitious seminar list at the department entrance cannot exist without coursework and student inquiry beneath it. Integrating the interpretation and verification of machine reasoning into this training does not compromise the tradition, but extends it, preparing students for a world in which the capacity to audit a mathematically sophisticated argument is scarcer, and therefore more valuable, than the capacity to generate one.

\smallskip

These proposals do not contradict the defense of ``useless'' exploratory mathematics; they depend on it. Auditing, comprehension, and structural translation cannot be taught in isolation as narrow procedural skills; they are byproducts of the deep, undirected inquiry that pure mathematics sustains, forged by contending with objective structures for their own sake. The assurance offices, residencies, and fellowships simply connect that discernment to the institutions that need it most, without asking the discipline to surrender its exploratory core.

\smallskip

We stand at the beginning of a profound transformation. The question is whether, when that transformation is complete, there will still be human beings who can look at a complex chain of mathematical reasoning, generated by any AI system, and say, with the authority and conceptual clarity earned through years of practice: \textit{this is what the result establishes, and what it does not, and here is why}.

\smallskip

If the answer is yes, the wandering continues; the borders hold, and our civilization keeps its sight.

\smallskip

If the answer is no, we will have all the theorems we could ever want, proved by machines of extraordinary power, verified by systems of impeccable precision. And no one will understand any of them. And no one, when the system fails, will know how to ask the right question.

\vspace{1.5em}
\begin{center}
*\quad*\quad*
\end{center}
\vspace{1em}

Somewhere tomorrow morning a graduate student will stop in front of a seminar list and read a title she does not yet understand. Everything this essay defends begins in that moment.

\newpage

\begin{center}
{\large\bfseries PURE MATHEMATICS TRUTHS}\\[0.5em]
{\footnotesize\itshape Appendix to \textrm{The Wandering Must Continue} --- Adapted from the principles formulated by John M.\ Collins~\cite{collins1987}}
\end{center}
\vspace{0.2em}
\hrule
\vspace{1.2em}

\bigskip

\noindent\textbf{Truth I: Minds are more important than compute.}

\noindent People --- not AI models --- make the critical difference. The right minds, rigorously trained and working together, will frame the right questions and craft the definitions that let us see the underlying structure. On the other hand, the most powerful compute in the world cannot compensate for a lack of humans who understand why a question matters or why a theorem holds.

\bigskip
\bigskip

\noindent\textbf{Truth II: Quality is better than quantity.}

\noindent A tight-knit group of mathematicians, brought together organically, deeply trained in their specialties, and guided by structural vision, is preferable to swarms of automated agents grinding through isolated conjectures without comprehension. A single unifying framework with a deep horizon outweighs a thousand opaque proofs.

\bigskip
\bigskip

\noindent\textbf{Truth III: Pure mathematicians cannot be mass produced.}

\noindent It takes years to train a human mind to reason in pure abstraction. Intense apprenticeship, in foundational graduate courses as well as active research seminars and conferences, is required to turn motivated students into working mathematicians. This process cannot be hastened without degrading the depth of the discipline.

\bigskip
\bigskip

\noindent\textbf{Truth IV: Mathematical capability cannot be created after emergencies occur.}

\noindent Building a living mathematical culture takes generations. The abstract structures that resolve a scientific crisis or lay the foundation for the next paradigm cannot be invented on demand; they must be discovered for their own sake, in active pure mathematics departments long before the need becomes acute.

\bigskip
\bigskip

\noindent\textbf{Truth V: Pure mathematics requires institutional support.}

\noindent The enduring strength of pure mathematics cannot be, and never has been, achieved without the broader ecosystem that sustains it. Modern mathematics is made possible by the undergraduate programs that draw in talented minds, the faculty who mentor them one by one, the neighboring sciences that pass questions back and forth with it, and the universities that shelter disinterested inquiry from short-term utilitarian pressure. Removing this support does not merely inconvenience the endeavor --- it makes the work impossible.

\newpage

\noindent\textbf{Acknowledgments.}
The ideas in this essay grew out of the author's research in arithmetic geometry and Hodge theory and cross-disciplinary collaboration in neuroscience. The author thanks Arun Ram for valuable suggestions on an earlier draft. Generative artificial intelligence tools (Claude, Anthropic; and Gemini, Google) were used for editorial feedback and prose refinement during the preparation of this manuscript. All arguments and conclusions are the author's own.

\vspace{+10pt}

\noindent Jun--Yong Park \enspace ---\enspace \texttt{june.park@sydney.edu.au} \\
\textsc{School of Mathematics and Statistics, University of Sydney, Australia}

\end{document}